\documentclass{amsart}
\usepackage[cp1251]{inputenc}
\usepackage[english,russian]{babel}
\usepackage{amsmath}
\usepackage{amssymb}
\usepackage{amsfonts}

\newtheorem{lem}{Лемма}
\newtheorem{thm}{Теорема}

\begin{document}
\thispagestyle{empty}

\title {Heat equation, Dirichlet problem, Laplace equation and randomized probability density
  }% Specify the name of the article
 \author { Oleg Yaremko,Alexander Abuzov, Dina Khotsyan}

\address{Oleg Yaremko,
\newline\hphantom {iii} Penza State Pedagogical University,% Place of work
\newline\hphantom {iii} st. Krasnaya, 40% address (street, house, structure, etc.)
\newline\hphantom {iii} 440038, Penza, Russia}%
\address{Alexander Abuzov,
\newline\hphantom {iii} Penza State Pedagogical University,% Place of work
\newline\hphantom {iii} st. Krasnaya, 40% address (street, house, structure, etc.)
\newline\hphantom {iii} 440038, Penza, Russia}%
\address{Dina Khotsyan,
\newline\hphantom {iii} Penza State Pedagogical University,% Place of work
\newline\hphantom {iii} st. Krasnaya, 40% address (street, house, structure, etc.)
\newline\hphantom {iii} 440038, Penza, Russia}%
\maketitle {\small
\begin{quote}

\noindent{\bf Abstract } In the present study examines the statistical structure of the generated randomized density of the normal distribution and the Cauchy distribution. The study put the allegation that a randomized probability density of the normal distribution can be regarded as the solution of the Cauchy problem for the heat equation, and randomized probability density of the Cauchy distribution can be considered as a solution to the Dirichlet problem for the Laplace equation. Conversely, the solution of the Cauchy problem for the heat equation can be regarded as a randomized probability density of the normal distribution, and the solution of the Dirichlet problem for the Laplace equation as randomized probability density of the Cauchy distribution. The main objective of the study was the fact that for each of these two cases to find the Fisher information matrix components and structural tensor.  We found nonlinear differential equations of the first, second and third order for the density of the normal distribution and Cauchy density computational difficulties to overcome . The components of the metric tensor (the Fisher information matrix) and the components of the strain tensor are calculated according to formulas in which there is the log-likelihood function, ie, logarithm of the density distribution. Because of the positive definiteness of the Fisher information matrix obtained inequality, which obviously satisfy the Cauchy problem solution with nonnegative initial conditions in the case of the Laplace equation and the heat equation.
\medskip

\noindent{\bf Key words:}Fisher information matrix,  structure tensor,  random density,Poisson formula,  Heat equation, Dirichlet problem, Laplace equation.
\end{quote} }
\section{introduction}

В настоящее время широко используются методы дифференциальной геометрии в исследовании информационных массивов (семейств вероятностных распределений пространств квантовых состояний, нейронных сетей и т.п.). Исследования по информационной геометрии восходят к статье С.Рао [1], где на основе фишеровской информационной матрицы была определена риманова метрика на многообразии распределений вероятностей. Дальнейшие исследования привели к понятию статистического многообразия как гладкого $n$ мерного многообразия на $M$ на котором задана метрически-аффинная структура $(g,\bar{\nabla})$, где $g$ -риманова метрика, а $\bar{\nabla}$-линейная связность без кручения., совместимая с метрикой $g$, т.е. для $g$ и $\bar{\nabla}$ выполняется условие Кодацци []
\[(\bar{\nabla}_{x}g(y,z)=\bar{\nabla}_{y}g(x,z)
\]
для любых векторных $x,y,z$ на $M$. Такая метрически-аффинная структура называется статистической. Так как $\nabla=\bar{\nabla}+\bar{T}$ , где $\nabla$ связность Леви-Чевита метрики $g$, а $\bar{T}$ ее тензор деформации , то ковариантный тензор деформации $T(x,y,z)=g((T(x,y),z)$, в силу условия Кодацци, симметричен по своим аргументам. Таким образом статистическая структура определяется заданием на $M$ пары тензорных полей $(g,T)$.

Статистической моделью называется гладко параметризованное конечным числом действительных параметров $\theta^{i}(i=\bar{1,n})$ семейство $S$распределений вероятностей $P_{\theta}|\theta\in R^{n}$ случайной величины. Всякое распределение вероятностей $P_{\theta}$ случайной величины $\xi$ характеризуется своей плотностью $p(\xi|\theta)$ на выборочном пространстве $\Omega,\xi \in \Omega$ или некоторой функцией от $p$ из непрерывного 1 -параметрического семейства функций
\[
\varphi_{2}(p)=\left( {\begin{array}{*{20}{c}}
   \frac{\alpha}{1-\alpha} p^{\frac{1-\alpha}{\alpha}},\alpha\neq1 \\
   { \ln p},\alpha=1
\end{array}} \right)
\]
При $\alpha=1$получаем функцию правдоподобия $\varphi=\ln p$. В этом случае компоненты  $g_{ij}$ метрического тензора $g$ (информационная матрица Фишера) и компоненты $T_{ijk}$ тензора деформации вычисляются по следующим формулам:
\begin{equation}
g_{ij}(\theta)=\int\partial_{i}\ln p\cdot\partial_{j}\ln p\cdot p d\xi
\end{equation}
\begin{equation}
T_{ijk}(\theta)=-\frac{1}{2}\int\partial_{i}\ln p\cdot\partial_{j}\ln p\cdot \partial_{k}\ln p\cdotp d\xi
\end{equation}
где $\partial_{i}=\partial / \partial\theta^{i}$.

\section{Рандомизированная  плотность нормального распределения и соответствующая статистическая структура }
Рассмотрим уравнение теплопроводности с начальным условием
\begin{equation}
\begin{cases}
u_t(x,t) -  u_{xx}(x,t) = 0& (x, t) \in \mathbf{R} \times (0, \infty)\\
u(x,0)=f(x)&
\end{cases}
\end{equation}
где $f(x)$ начальное распределение температурного поля. Решение задачи Коши выражается через фундаментальное решение
\[
\Phi(x,t)=\frac{1}{2\sqrt{\pi t}}\exp\left(-\frac{x^2}{4t}\right)
\]
при помощи следующей формулы:
\begin{equation}
u(x,t) = \int \Phi(x-\xi,t) f(\xi) d\xi.
\end{equation}

Замечание. Если $f(\xi)$-некоторая плотность вероятности, то решение задачи Коши (1) можно интерпретировать как рандомизированную  плотность [] нормального распределения с параметрами $(x,\sqrt{2t})$.

Обозначим $$h(t,x-\xi)=2\sqrt{\pi t}\Phi(x-\xi,t)=e^{-\frac{(x-\xi)^{2}}{4t}}.$$
\begin{lem} Справедливы соотношения $$\frac{h'_{t}h'_{t}}{h}=h''_{tt}+\frac{2}{t}h'_{t}, \frac{h'_{x}h'_{x}}{h}=h''_{xx}+\frac{1}{2t}h, \frac{h'_{t}h'_{x}}{h}=h''_{tx}+\frac{1}{t}h'_{x}.$$
\end{lem}
Доказательство.
$$h'_{x}=-\frac{x-\xi}{2t}h,h'_{x}h'_{x}=\frac{(x-\xi)^{2}}{4t^{2}}h^{2}.$$
Тогда $$h''_{xx}=-\frac{1}{2t}h+\frac{(x-\xi)^{2}}{4t^{2}}h,
\frac{h'_{x}h'_{x}}{h}=h''_{xx}+\frac{1}{2t}h,$$
$$h'_{t}=\frac{(x-\xi)^{2}}{4t^{2}}h,
h'_{t}h'_{t}=\left(\frac{(x-\xi)^{2}}{4t^{2}}\right)^{2}h^{2},$$
$$h''_{tt}=-\frac{(x-\xi)^{2}}{2t^{3}}h+\left(\frac{(x-\xi)^{2}}{4t^{2}}\right)^{2}h.$$
Значит,
$$\frac{h'_{t}h'_{t}}{h}=h''_{tt}+\frac{2}{t}h'_{t}.$$
Аналогично,
$$\frac{h'_{t}h'_{x}}{h}=h''_{tx}+\frac{1}{t}h'_{x}.$$
\begin{lem} Выполняются следующие тождества
$$\frac{h'^{3}_{x}}{h^{2}}=h'''_{xxx}+\frac{3}{2t}h'_{x}$$
$$\frac{h'^{2}_{x}h'_{t}}{h^{2}}=h'''_{xxt}+\frac{2}{t}h''_{xx}+\frac{1}{2t}h'_{t}+\frac{1}{2t^{2}}h$$
$$\frac{h'^{3}_{t}}{h^{2}}=h'''_{ttt}+\frac{6}{t}h''_{tt}+\frac{6}{t^{2}}h'_{t}$$
$$\frac{h'^{2}_{t}h'_{x}}{h^{2}}=h'''_{xtt}+\frac{4}{t}h''_{tx}+\frac{2}{t^{2}}h'_{x}$$
\end{lem}
Доказательство. Достаточно продифференцировать по переменной $t$ или $x$ каждое из равенств в лемме 1.

\begin{lem} Пусть функция $f(x)$ неотрицательна на действительной оси, а функция $u(t,x)$ определяется формулой (1), тогда функция
$$p(\xi,x,t)=\frac{h(t,x-\xi)f(\xi)}{2\sqrt{\pi t}u(t,x)}$$
является плотность распределения,  здесь$\theta^{1}=x,\theta^{2}=t$- параметры семейства $S$ распределений вероятностей, $\xi$ -случайная величина на выборочном пространстве $R$.
\end{lem}
Найдем статистическую структуру на $S$ в случае, когда $\varphi$ является функцией правдоподобия, т.е. $\varphi=\ln p$. В этом случае информационная матрица Фишера вычисляется явно.

\begin{thm} Элементы информационной матрицы Фишера-компоненты метрического тензора $g$ имеют вид
$$g_{11}=(\ln u)''_{xx}+\frac{1}{2t}$$
$$g_{22}=(\ln u)''_{tt}+\frac{2}{t}(\ln u)'_{t}+\frac{1}{2t^{2}}$$
$$g_{12}=(\ln u)''_{xt}+\frac{1}{t}(\ln u)'_{x}$$
\end{thm}
Доказательство. Вычислим информационную матрицу
$$g_{11}=\int (\ln p)'_{x}(\ln p)'_{x}pd\xi=\int \left((\ln h)'_{x}(\ln h)'_{x}-2(\ln h)'_{x}(\ln u)'_{x}+(\ln u)'_{x}(\ln u)'_{x}\right)\frac{hf}{u}=$$
$$=\int\frac{h'_{x}h'_{x}}{h}\frac{f}{u}d\xi-2(\ln u)'_{x}\int\frac{h'_{x}f}{u}d\xi+(\ln u)'_{x}(\ln u)'_{x}$$

По формулам () вычисляем слагаемые
$$\int\frac{h'_{x}f}{u}d\xi=\frac{u'_{x}}{u}$$
$$\int\frac{h'_{x}h'_{x}}{h}\frac{f}{u}d\xi=\int \left(h''_{xx}+\frac{1}{2t}h \right)\frac{f}{u}d\xi=\frac{u''_{xx}}{u}+\frac{1}{2t}$$
В результате
$$g_{11}=(\ln u)''_{xx}+\frac{1}{2t}$$
Аналогично
$$g_{12}=\int (\ln p)'_{x}(\ln p)'_{t}pd\xi=\int ((\ln h)'_{x}(\ln h)'_{t}-(\ln h)'_{x}(\ln \sqrt{t}u)'_{t}-$$
$$
-(\ln h)'_{t}(\ln \sqrt{t}u)'_{x}+(\ln \sqrt{t}u)'_{x}(\ln \sqrt{t}u)'_{t})\frac{hf}{\sqrt{t}u}=$$
$$=\int\frac{h'_{x}h'_{t}}{h}\frac{f}{\sqrt{t}u}d\xi-(\ln \sqrt{t}u)'_{t}\int\frac{h'_{x}f}{\sqrt{t}u}d\xi-(\ln \sqrt{t}u)'_{x}\int\frac{h'_{t}f}{\sqrt{t}u}d\xi+(\ln \sqrt{t}u)'_{x}(\ln \sqrt{t}u)'_{t}=
$$
$$=(\ln \sqrt{t}u)''_{xt}+\frac{1}{t}(\ln \sqrt{t}u)'_{x}=(\ln u)''_{xt}+\frac{1}{t}(\ln u)'_{x}$$
Точно так же
$$g_{22}=\int (\ln p)'_{t}(\ln p)'_{t}pd\xi=\int \left((\ln h)'_{t}(\ln h)'_{t}-2(\ln h)'_{t}(\ln \sqrt{t}u)'_{t}+(\ln \sqrt{t}u)'_{t}(\ln \sqrt{t}u)'_{t}\right)\frac{hf}{\sqrt{t}u}=$$
$$=\int\frac{h'_{t}h'_{t}}{h}\frac{f}{\sqrt{t}u}d\xi-2(\ln \sqrt{t}u)'_{t}\int\frac{h'_{t}f}{\sqrt{t}u}d\xi+(\ln \sqrt{t}u)'_{t}(\ln \sqrt{t}u)'_{t}=$$
$$=(\ln  \sqrt{t}u)''_{tt}+\frac{2}{t}(\ln \sqrt{t}u)'_{t}=(\ln u)''_{tt}+\frac{2}{t}(\ln u)'_{t}+\frac{1}{2t^{2}}$$
\begin{thm}Компоненты структурного тензора имеют вид
$$T_{111}=-\frac{1}{2}(\ln u)'''_{xxx}-\frac{3}{4t}(\ln u)'_{x}$$
$$T_{112}=-\frac{1}{2}(\ln u)'''_{xxt}-\frac{1}{t}(\ln u)''_{xx}-\frac{1}{4t}(\ln u)'_{t}-\frac{1}{4t^{2}}$$
$$T_{122}=-\frac{1}{2}(\ln u)''_{xtt}-\frac{2}{t}(\ln u)''_{tx}-\frac{1}{t^{2}}(\ln u)'_{x}$$
$$T_{222}=-\frac{1}{2}(\ln u)'''_{ttt}-\frac{3}{t}(\ln u)''_{tt}-\frac{3}{t^{2}}(\ln u)'_{t}$$
\end{thm}
Доказательство.
\[{T}_{222}(\theta)=-\frac{1}{2}\int
 \partial_{t} \log p(\xi;t,x)\cdot \partial_{t}\log p(\xi;t,x)\,\cdot\partial_{t}\log p(\xi;t,x)\cdot p(\xi;t,x)d\xi=
\]

\[=-\frac{1}{2}\int\left(\frac{h'_{t}}{h}-\frac{v'_{t}}{v}\right)^{3}\frac{hf(\xi)}{v}d\xi=
\]
\[=-\frac{1}{2}\int\left(\frac{h'^{3}_{t}}{h^{3}}-3\frac{h'^{2}_{t}}{h^{2}}\frac{v'_{t}}{v}+3\frac{h'_{t}}{h}\frac{v'^{2}_{t}}{v^{2}}-\frac{v'^{3}_{t}}{v^{3}}\right)\frac{hf(\xi)}{v}d\xi
\]
Применим лемму 2 и воспользуемся равенствами
\[\int h(t,x-\xi)f(\xi)d\xi=v(t,x),\int h'_{t}(t,x-\xi)f(\xi)d\xi=v'_{t}(t,x),
\]
\[\int h''_{tt}(t,x-\xi)f(\xi)d\xi=v''_{tt}(t,x),\int h'''_{ttt}(t,x-\xi)f(\xi)d\xi=v'''_{ttt}(t,x),v=2\sqrt{\pi t}u\]
Элементарными преобразованиями получим требуемую компоненту структурного тензора. Остальные компоненты вычисляются аналогично.

\section{Рандомизированная  плотность  распределения Коши и статистическая структура}
Рассмотрим задачу Дирихле для уравнения Лапласа для полуплоскости
\begin{equation}
\begin{cases}
u_{xx}(x,y) + u_{yy}(x,y) = 0& (x, y) \in  (0, \infty)\times\mathbf{R} \\
u(0,y)=f(y)&
\end{cases}
\end{equation}
где $f(y)$ граничное значение. Решение задачи Дирихле выражается через фундаментальное решение
\[
\Phi(x,y)=\frac{1}{\pi}\frac{x}{x^{2}+y^{2}}
\]
при помощи следующей формулы:
\begin{equation}
u(x,t) = \int \Phi(x,y-\xi) f(\xi) d\xi.
\end{equation}

Замечание. Если $f(\xi)$-некоторая плотность вероятности, то решение задачи Дирихле (3) можно интерпретировать как рандомизированную  плотность [] распределения Коши.
 Обозначим $$h(x,y-\xi)=\frac{\pi \Phi(x,y-\xi)}{x}=\frac{1}{(x^{2}+(y-\eta)^{2})},$$
\begin{lem} Плотность распределения Коши удовлетворяет нелинейным дифференциальным уравнениям второго порядка
$$\frac{h'_{x}h'_{x}}{h}=h''_{xx}-\frac{1}{x}h'_{x}$$
$$\frac{h'_{y}h'_{y}}{h}=h''_{yy}-\frac{1}{x}h'_{x}$$
$$\frac{h'_{y}h'_{x}}{h}=\frac{1}{2}h''_{yx}$$
\end{lem}
Доказательство.
Имеем соотношения $$h'_{x}=-\frac{2x}{x^{2}+(y-\eta)^{2}}h,h'_{x}h'_{x}=\frac{4x^{2}}{(x^{2}+(y-\eta)^{2})^{2}}h^{2}$$
Тогда $$h''_{xx}=-\frac{2}{x^{2}+(y-\eta)^{2}}h+\frac{4x^{2}}{(x^{2}+(y-\eta)^{2})^{2}}h$$
$$\frac{h'_{x}h'_{x}}{h}=h''_{xx}-\frac{1}{x}h'_{x}$$
$$h'_{y}=-\frac{2(y-\eta)}{x^{2}+(y-\eta)^{2}}h,h'_{y}h'_{y}=\frac{4(y-\eta)^{2}}{(x^{2}+(y-\eta)^{2})^{2}}h^{2}$$

$$h''_{yy}=-\frac{2}{x^{2}+(y-\eta)^{2}}h+\frac{4(y-\eta)^{2}}{(x^{2}+(y-\eta)^{2})^{2}}h$$
Значит
$$\frac{h'_{y}h'_{y}}{h}=h''_{yy}-\frac{1}{x}h'_{x}$$
Аналогично
$$\frac{h'_{y}h'_{x}}{h}=\frac{1}{2}h''_{yx}$$
\begin{lem} Плотность распределения Коши удовлетворяет нелинейным дифференциальным уравнениям третьего порядка
$$\frac{h'^{3}_{x}}{h^{2}}=h'''_{xxx}-\frac{3}{x}h''_{xx}+\frac{3}{x^{2}}h'_{x}$$
$$\frac{h'^{3}_{y}}{h^{2}}=h'''_{yyy}-\frac{3}{2x}h''_{xy}$$
$$\frac{h'^{2}_{y}h'_{x}}{h^{2}}=\frac{1}{3}h'''_{yyx}-\frac{1}{3x}h''_{xy}+\frac{1}{3x^{2}}h'_{x}$$
$$\frac{h'^{2}_{x}h'_{y}}{h^{2}}=\frac{1}{3}h'''_{xxy}-\frac{1}{3x}h''_{xy}$$
\end{lem}
\begin{lem} Пусть функция $f(x)$ неотрицательна на действительной оси. тогда формула
$$p(\xi,x,y)=\frac{x h(t,x-\xi)f(\xi)}{\pi u(x,y)}$$
определяет плотность распределения.
\end{lem}
\begin{thm} Элементы информационной матрицы имеют вид
$$g_{11}=(\ln u)''_{xx}-\frac{1}{x}(\ln u)'_{x}+\frac{2}{x^{2}}$$
$$g_{12}=\frac{1}{2}(\ln u)''_{xy}$$
$$g_{22}=(\ln u)''_{xx}-\frac{1}{x}(\ln u)'_{x}+\frac{2}{x^{2}}$$
\end{thm}
Доказательство.
$$g_{11}=\int (\ln p)'_{x}(\ln p)'_{x}pd\xi=\int \left((\ln h)'_{x}(\ln h)'_{x}-2(\ln h)'_{x}(\ln v)'_{x}+(\ln v)'_{x}(\ln v)'_{x}\right)\frac{hf}{v}=$$
$$=\int\frac{h'_{x}h'_{x}}{h}\frac{f}{v}d\xi-2(\ln v)'_{x}\int\frac{h'_{x}f}{v}d\xi+(\ln v)'_{x}(\ln v)'_{x},$$
где $v=\frac{\pi}{x}u.$
По лемме вычисляем слагаемые
$$\int\frac{h'_{x}f}{v}d\xi=\frac{v'_{x}}{v}$$
$$\int\frac{h'_{x}h'_{x}}{h}\frac{f}{v}d\xi=\int \left(h''_{xx}-\frac{1}{x}h'_{x} \right)\frac{f}{v}d\xi=\frac{v''_{xx}}{v}-\frac{1}{x}\frac{v'_{x}}{v}$$
В результате находим
$$g_{11}=(\ln v)''_{xx}-\frac{1}{x}(\ln v)'_{x}$$
Аналогично
$$g_{12}=\int (\ln p)'_{x}(\ln p)'_{y}pd\xi=\int \left((\ln h)'_{x}(\ln h)'_{y}-(\ln h)'_{x}(\ln v)'_{y}-(\ln h)'_{t}(\ln v)'_{x}+(\ln v)'_{x}(\ln v)'_{y}\right)\frac{hf}{v}=$$
$$=\int\frac{h'_{x}h'_{y}}{h}\frac{f}{v}d\xi-(\ln v)'_{y}\int\frac{h'_{x}f}{v}d\xi-(\ln v)'_{x}\int\frac{h'_{y}f}{v}d\xi+(\ln v)'_{x}(\ln v)'_{y}=\frac{1}{2}\frac{v''_{xy}}{v}-\frac{v'_{x}v'_{y}}{v^{2}}$$
$$g_{22}=\int (\ln p)'_{y}(\ln p)'_{y}pd\xi=\int \left((\ln h)'_{y}(\ln h)'_{y}-2(\ln h)'_{y}(\ln v)'_{y}+(\ln v)'_{y}(\ln v)'_{y}\right)\frac{hf}{v}=$$
$$=\int\frac{h'_{y}h'_{y}}{h}\frac{f}{v}d\xi-2(\ln v)'_{y}\int\frac{h'_{y}f}{v}d\xi+(\ln v)'_{y}(\ln v)'_{y}=$$
$$=(\ln v)''_{xx}-\frac{1}{x}(\ln v)'_{x}$$

\begin{thm} Компоненты структурного тензора имеют вид
$$T_{111}=-\frac{1}{2}(\ln u)'''_{xxx}+\frac{3}{2x}(\ln u)''_{xx}-\frac{3}{2x^{2}}(\ln u)'_{x}-\frac{4}{x^{3}}$$
$$T_{112}=-\frac{1}{6}(\ln u)'''_{xxy}+\frac{1}{6x}(\ln u)''_{xy}$$
$$T_{122}=-\frac{1}{6}(\ln u)'''_{yyx}+\frac{1}{6x}(\ln u)''_{xy}-\frac{1}{6x^{2}}(\ln u)'_{x}+\frac{1}{6x^{3}}$$
$$T_{222}=-\frac{1}{2}(\ln u)'''_{yyy}+\frac{3}{4x}(\ln u)''_{xy}$$
\end{thm}
Рассуждения проводятся по образцу теоремы 1 с использованием лемм 4,5,6.
\section{Заключение}
Нами проведено исследование  статистических структур, порождаемых рандомизированными  плотностями нормального распределения и распределения Коши. В основу исследования положено утверждение о том,  что рандомизированную  плотность вероятности нормального распределения можно рассматривать как решение задачи Коши для уравнения теплопроводности, а  рандомизированную  плотность вероятности  распределения Коши можно рассматривать как решение задачи Дирихле для  уравнения Лапласа. Обратно, решение задачи Коши для уравнения теплопроводности можно рассматривать как рандомизированную  плотность вероятности нормального распределения, а  решение задачи Дирихле для  уравнения Лапласа как рандомизированную  плотность вероятности  распределения Коши. Для каждого из этих двух случаев мы нашли компоненты  информационной матрицы Фишера и структурного тензора. Предложен новый метод вычисления этих компонент, основанный на выведенных нами нелинейных дифференциальных уравнениях первого, второго и третьего порядков для плотностей нормального распределения и плотности Коши. В качестве следствия из положительной определенности информационной матрицы Фишера, можно получить неравенства, которым заведомо удовлетворяют решения задачи Коши для уравнения теплопроводности с неотрицательным начальным условием и решения задачи Дирихле для уравнения Лапласа с неотрицательным краевым значением  для случая уравнения Лапласа.

  \end{document}